\documentclass[11pt,epsf]{article}
\usepackage{ifpdf}
\usepackage{amsthm}
\usepackage{amsmath}
\usepackage{amssymb}
\usepackage{graphicx}
\usepackage{psfrag}
\usepackage[all]{xypic}
\usepackage{url}

\newcommand{\sinpioverthree}{0.866025404}

\newcommand{\frkF}{\mathfrak{F}}
\newcommand{\frkG}{\mathfrak{G}}
\newcommand{\tand}{{\tiny \text{AND} }}
\newcommand{\nbhd}{\operatorname{nbhd}}
\newcommand{\disjointunion}{\sqcup}
\newcommand{\isfaceof}{\prec}

\newcommand{\R}{\mathbb{R}}

\newcommand{\N}{\mathbb{N}}

\newcommand{\union}{\cup}

\newcommand{\indep}{\! \perp \!\!\! \perp \!}

\newtheorem{theorem}{Theorem}
\newtheorem{lemma}[theorem]{Lemma}

\newtheorem{proposition}[theorem]{Proposition}
\newtheorem{observation}[theorem]{Observation}

\theoremstyle{definition}
\newtheorem{definition}[theorem]{Definition}
\newtheorem{example}[theorem]{Example}
\newtheorem{remark}[theorem]{Remark}
\newtheorem{algorithm}[theorem]{Algorithm}

\title{Convex Rank Tests and Semigraphoids}
\date{}
\author{Jason Morton, Lior Pachter, Anne Shiu, Bernd Sturmfels, and Oliver Wienand}
\begin{document}

\maketitle

\begin{abstract}
Convex rank tests are partitions of the symmetric group
which have desirable geometric properties.
The statistical tests defined by such partitions involve counting all
permutations in the equivalence classes.
Each class consists of the linear extensions of a partially
ordered set specified by data.
Our methods refine existing rank tests of non-parametric statistics, such as
the sign test and the runs test, and are useful
for exploratory analysis of ordinal data.
We establish a bijection between  convex rank tests and probabilistic
conditional independence structures known as semigraphoids. The subclass of
submodular rank tests is derived from faces of the cone
of submodular functions, or from Minkowski summands of the permutohedron.
We enumerate all small instances of such rank tests. Of particular interest 
are  graphical tests, which correspond to both graphical models and to graph 
associahedra.
\\[2ex]
\end{abstract}

\section{Introduction}

The non-parametric approach to statistics was introduced by
\cite{Pitman1937SignificanceI} 
via the method of permutation testing.
Subsequent development of these ideas revealed a close connection
between non-parametric tests and {\em rank tests}, which are
statistical tests suitable for ordinal data.
Beginning in the 1950s,  many rank tests were developed for specific
applications, such as the comparison of populations or testing hypotheses for determining
the location of a population. The geometry of these tests was explored in \cite{Cook}. More recently, the search for patterns
in large datasets has spurred the development and exploration of new
tests. For instance, the emergence of microarray data in
molecular biology has led to tests for identifying significant
patterns in gene expression time series; see e.g.~\cite{Willbrand2005}.
This application motivated us to
develop a mathematical theory of rank tests.
We propose
that a {\em rank test} is a partition of $S_n$ induced by a
map $\, \tau : S_n \rightarrow T\,$
from the symmetric group of all permutations of $[n]=\{1,\ldots,n\}$
onto a set $T$ of statistics.
The statistic $\tau(\pi)$ is the {\em signature} of the permutation $\pi \in S_n$.
Each rank test defines a partition of $S_n$ into
classes, where $\pi$ and $\pi'$ are in the same class if and only if
$\tau (\pi) = \tau(\pi')$. We identify $T = {\rm image}(\tau)$ with
the set of all classes in this partition of $S_n$.
 Assuming the uniform distribution on  $ S_n$, the probability
of seeing a particular signature $t \in T$ is
$\,1/n! \,$ times $| \tau^{-1}( t)|$.
The computation of a $p$-value for a given permutation $\pi \in S_n$
leads to the problem of summing
\begin{equation} \label{Pvalue}
{\rm Pr}(\pi') \quad = \quad
 \frac{1}{n !} \cdot |\, \tau^{-1} \bigl( \tau(\pi') \bigr)\, |
 \end{equation}
over permutations $\pi'$ with ${\rm Pr}(\pi') \leq {\rm Pr}(\pi)$,
a computational task to be addressed in Section~6.

This paper is an expanded version of our note ``Geometry of Rank Tests''
which was presented in September 2006 in Prague at the conference
{\em Probabilistic Graphical Models (PGM 3)}.
The emphasis of our discussion is on the mathematics underlying rank tests, and, in particular,
on the connection to statistical learning theory (semigraphoids).
We refer to \cite{cyclohedron} for  details on how to use our rank tests in practice, and how to
interpret the p-values derived from (\ref{Pvalue}).

The five subsequent sections are organized as follows.
 In Section 2 we explain how existing rank tests in
 non-parametric statistics  can be understood from our
 geometric  point of view, and how they are described
 in the language of algebraic combinatorics \cite{Stanley1997}.
In Section 3  we define the class of {\em  convex rank tests}.
These tests are most natural from both the statistical and
the combinatorial point of view.  Convex rank tests can
be defined as polyhedral fans that coarsen the 
hyperplane arrangement of $S_n$.
Our main result (Theorem \ref{fantheorem})
states that convex rank tests are in bijection with
 conditional independence structures known as
{\em semigraphoids} \cite{Dawid, Pearl, Studeny2005Probabilistic}.

Section 4 is devoted to convex rank tests that are
induced by submodular functions.
These {\em submodular rank tests} are in
bijection with Minkowski summands
of the $(n{-}1)$-dimensional permutohedron and with {\em structural
imset models}. These tests are
at a suitable level of generality for the biological applications
\cite{cyclohedron, Willbrand2005} that motivated us.
The connection between polytopes and independence models is
made concrete in the classification of small models in Remarks
\ref{rmk1}--\ref{rmk3}.

In Section 5 we study the subclass of {\em graphical tests}.
In combinatorics, these correspond to
graph associahedra, and in statistics
to graphical models.
The equivalence of these two structures is shown in
Theorem \ref{maingraphical}.
The implementation of convex rank tests requires the efficient enumeration of
linear extensions of partially ordered sets.  
Our algorithms and software are discussed in Section~6.
A key ingredient is the efficient computation of distributive lattices.

\section*{Acknowledgments}
Our research on rank tests originated in discussions with Olivier Pourqui\'{e} and Mary-Lee Dequ\'{e}ant as part of
the DARPA Program {\em Fundamental Laws of Biology}, that supported Jason Morton, Lior Pachter, and Bernd Sturmfels.  Anne Shiu was supported by a Lucent Technologies Bell Labs Graduate Research Fellowship, and Oliver Wienand by the Wipprecht foundation.  We thank Milan Studen\'{y} and Franti\v{s}ek Mat\'{u}\v{s} for helpful comments.

\bigskip

\section{Rank tests and posets}

A permutation $\pi$ in $S_n$ is a
total order on the set $[n] := \{1,\ldots,n\}$.
This means that $\pi$ is a set
of $\binom{n}{2}$ ordered pairs
of elements in $[n]$.  For example, $\pi = \{ (1,2), (2,3), (1,3) \}$ represents the total order $1>2>3$.
If $\pi$ and $\pi'$ are permutations then
$\,\pi \cap \pi'\,$ is a partial order.

In the applications we have in mind, the data
are vectors $u \in \R^n$ with distinct coordinates.
The permutation associated with $u$ is the
total order $\,\pi = \{ \,(i,j)\in [n] \times [n] \,: \, u_i < u_j\,\}$.
We shall employ two other ways of writing a permutation.
The first is the {\em rank vector} $\,\rho = (\rho_1,\ldots,\rho_n)$,
whose defining properties are
$\{\rho_1,\ldots,\rho_n\} = [n]$ and
$\rho_i < \rho_j$ if and only if $u_i < u_j$.  That is, the coordinate of the rank vector with value $i$ is at the same position as the $i$th smallest coordinate of $u$.  
The second is the {\em descent vector}
$\delta = (\delta_1 | \delta_2 | \ldots | \delta_n)$.
The descent vector is defined by $u_{\delta_i} > u_{\delta_{i+1}}$
for $i=1,2,\ldots, n {-} 1$.  Thus the $i$th coordinate of the descent vector is the position of the $i$th largest value of the data vector $u$.
For example, if $\,u = (11,7,13)\,$ then its permutation
is represented by
$\, \pi = \{ (2,1),(1,3),(2,3)\}$, by $\,\rho = (2,1,3)$, or by  $\,\delta = (3|1|2)$.

A permutation $\pi $ is a {\em linear extension} of a
partial order $P$ on $[n]$ if $P \subseteq \pi$, i.e. $\pi$ is a total order that refines the partial order $P$.
 We write $\mathcal{L}(P) \subseteq S_n$ for the set of
linear extensions of~$P$.
A partition $\tau$ of the symmetric group $S_n$ is
a {\em pre-convex rank test}
if the following axiom holds:
\[\begin{array}{ccc}
{\rm (PC)} && \begin{array}{c} \text{If }\tau(\pi) = \tau(\pi') \text{ and } \pi'' \in \mathcal{L} (\pi \cap \pi')  \text{ then } \tau(\pi) \! = \! \tau(\pi') \!= \!\tau(\pi''). \end{array}\end{array}
\]
Note that $\, \pi'' \in \mathcal{L} (\pi \cap \pi') \,\,$  means
$\pi \cap \pi' \subseteq \pi''$.
The number of all rank tests $\tau$ on $[n]$ is the {\em Bell number} $B_{n!}$, which is
the number of set partitions of a set of cardinality $n!$. 

\begin{example} \label{bell}
For $n=3$ there are $B_6 = 203$ rank tests, or
partitions of the symmetric group $S_3$, which consists
of six permutations.
Of these $203$ rank tests, only $40$ satisfy the
axiom (PC).
One example is the pre-convex rank test in
Figure 1. Here
the symmetric group $S_3$  is partitioned
into the four classes
 $\,\bigl\{  (1|2|3)\bigr\}$,  $\,\bigl\{(2|1|3)\bigr\}$,
$\,\bigl\{(2|3|1)\bigr\}$, and
$\,\bigl\{(1|3|2), (3|1|2), (3|2|1)\bigr\}$.  
\end{example}

Each class $C$ of a pre-convex rank test $\tau$ corresponds to
a poset $P$ on the ground set $[n]$; namely, the partial order $P$ is the
intersection of all total orders in that class: $P=\bigcap_{\pi \in C} \pi$. The axiom
(PC) ensures that $C$ coincides with
the set $\mathcal{L}(P)$ of all linear extensions of $P$.  The inclusion $C \subseteq \mathcal{L}(P)$ is clear.  The proof of the reverse inclusion $\mathcal{L}(P) \subseteq C$ is based on the fact that,
from any permutation $\pi$ in $ \mathcal{L}(P)$, we can obtain any other $\pi'$ in $ \mathcal{L}(P)$ by a sequence of reversals $(a,b) \mapsto (b,a)$,
where each intermediate $\hat{\pi}$ is also in $ \mathcal{L}(P)$. 
Consider any $\pi_0 \in \mathcal{L}(P)$ and 
suppose that $\pi_1 \in  C$ differs by only one reversal  $(a,b)\in \pi_0$, $(b,a) \in \pi_1$.  Then
$(b,a) \notin P$, so there is some $\pi_2 \in C$
such that $(a,b) \in \pi_2$; thus, $\pi_0\in  \mathcal{L} (\pi_1 \cap \pi_2)$ by (PC).
This shows $\pi_0 \in C$.

A pre-convex rank test therefore can be characterized by 
an unordered collection of posets $P_1,P_2,\ldots,P_k$ on $[n]$ that
satisfies the property that  the symmetric group $S_n$ is the disjoint union of the subsets
$\mathcal{L}(P_1),\mathcal{L}(P_2), \ldots, \mathcal{L}(P_k)$.
This structure was discovered independently and studied
by Postnikov, Reiner and Williams \cite[\S 3]{PRW}
who used the term {\em complete fan of posets} for
what we shall call a convex rank test in Section 3.
The posets $\,P_1, P_2, \ldots, P_k\,$ that represent the classes in
a pre-convex rank test capture the
shapes of data vectors.  In graphical rank tests (Section \ref{sec:graphical}), this shape can be interpreted as a smoothed topographic map of the data vector.

\begin{example}[The sign test for paired data] \label{sign_test}
The \emph{sign test} is performed on data that are paired as
two vectors $u=(u_1,u_2, \dots,u_m)$ and
$ v = ( v_1, v_2, \dots,  v_m)$.  The null hypothesis
is that the median of the differences $u_i - v_i$ is 0.
The test statistic is the number of differences
that are positive. This test is a rank test, because
$u$ and $v$ can be transformed into the
overall ranks of the $n=2m$ values, and the rank vector
entries can then be compared. This test coarsens the convex rank test which is
the MSS test of Section 4 with $\,\mathcal{K} \,= \,\{\{1,m+1\},\{2, m+2\}, \dots \}$.
\end{example}

\begin{example}[Runs tests]
A \emph{runs test} can be used when there is a natural ordering on the data
points, such as in a time series. The data are
transformed into a sequence of `pluses' and `minuses,' and
the null hypothesis is that the number
of observed runs is
no more than that expected by chance.  
 Common types of runs tests include the
sequential runs test (`plus' if consecutive data points increase, `minus'  if they
decrease), and the runs test to check randomness of residuals, i.e.
deviation from a curve fit to the data.
A runs test is a coarsening of a convex rank test,
known as {\em up-down analysis}
\cite[\S 6.1.1]{Willbrand2005}, which is described
 in Example \ref{ex.updwn} below.
\end{example}

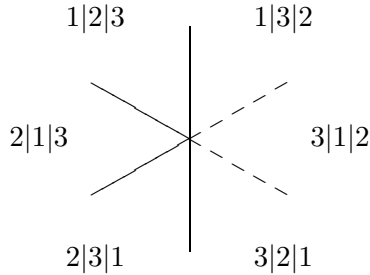
\begin{figure}[htb] \label{fig:nonconvexpreconvex}
\[
\begin{xy}<15mm,0mm>:
(0,0)  ="origin" ;
(\sinpioverthree, 0.5)  ="uprt"  ;
(\sinpioverthree, -0.5)  ="dwnrt"  ;
(-\sinpioverthree, 0.5)  ="upl" ;
(-\sinpioverthree, -0.5)  ="dwnl"  ;
(0, 1)  ="up" ;
(0,-1)  ="dwn" ;
   "origin";"uprt" **@{--};
   "origin";"dwnrt" **@{--};
   "origin";"upl" **@{-};
   "origin";"dwnl" **@{-};
   "origin";"up" **@{-};
   "origin";"dwn" **@{-};
(-.5,\sinpioverthree)  ="123"  *+!DR{1|2|3} ;
(.5 ,\sinpioverthree)  ="132"  *+!DL{1|3|2} ;
(1,0)                  ="312"  *+!L{3|1|2}  ;
(.5,-\sinpioverthree)  ="321"  *+!UL{3|2|1} ;
(-.5,-\sinpioverthree) ="231"  *+!UR{2|3|1} ;
(-1,0)                 ="213"  *+!R{2|1|3}  ;
 \end{xy}\]
\caption{Illustration of a pre-convex rank test that is not convex.  Cones are labelled by descent vectors, so $1|2|3$ indicates the cone $u_1 > u_2 > u_3$.  This rank test is specified by the
four posets $\,P_1 = \{3 {<} 1,2 {<} 1, 3 {<} 2\}, \, P_2 = \{1 {<}2, 3 {<} 2, 3 {<} 1\}, \,
 P_3 = \{3 {<}2, 1{<}3, 1{<}2\} \,$ and $\, P_4 = \{2 {<} 3\} $.}
\end{figure}

These two examples suggest
that many rank tests from
classical  non-parametric statistics  have a natural refinement by a
pre-convex rank test.  
However, not all tests have this property. 
Because many classical rank tests apply to loosely grouped data 
(e.g. data which are divided into two samples),  the
axiom (PC) is not always satisfied.
In such cases, the pre-convex rank test is a first step, after which 
permutations are grouped together under additional symmetries, e.g.,
the permutations  $\,\delta=\, (1|2|3|4|5)\,$ and 
 $\,\delta'=\, (5|4|3|2|1)\,$ might be identified.

The adjective ``pre-convex'' refers to the following interpretation
of the axiom (PC). Consider any two data vectors $u$ and $u'$
in $\R^n$, and a convex combination $u'' = \lambda u +
(1-\lambda) u'$, with $0 < \lambda < 1$.
If  $\,\pi, \pi', \pi'' \,$ are the permutations of $\, u,u',u'' \,$
then $\, \pi'' \in \mathcal{L}(\pi \cap \pi')$.
Thus the equivalence classes in $\R^n$ specified by a
 pre-convex rank test are convex cones.
 In the next section,  we shall remove the prefix 
 from ``pre-convex'' if the faces of these cones 
 fit together well.

\section{Convex rank tests}

A {\em fan} in $\R^n$ is a finite collection $\mathcal{F}$ of
polyhedral cones \cite{Ziegler1995} which satisfies the following properties:
\begin{itemize}
\item[(i)] if $C \in \mathcal{F}$ and $C'$ is a face of $C$, then
$C' \in \mathcal{F}$,
\item[(ii)] if $C, C' \in \mathcal{F}$,
then $C \cap C'$ is a face of $C$.
\end{itemize}
Two vectors $u$ and $v $ in $\R^n$ are
{\em permutation equivalent} when $u_i < u_j$
if and only if $v_i < v_j$, and $u_i = u_j$
if and only if $v_i = v_j$  for all $i,j \in [n]$.  Note that for two data vectors, each with distinct coordinates, they are permutation equivalent if and only if they have the same rank vector.
The permutation equivalence classes (of which there are $13$ for $n=3$) induce
a fan called the {\em $S_n$-fan}.  The arrangement of hyperplanes $\{x_i =x_j \}$ that defines these classes is also known as the {\em braid arrangement}, and its regions as the {\em Weyl chambers} of 
the Lie algebra $\mathfrak{sl}(n)$.  
The maximal cones in the $S_n$-fan, which are the closures
of the permutation equivalence classes,
are indexed
by permutations $\delta$ in $ S_n$.
A {\em coarsening} of the $S_n$-fan is a fan $\mathcal{F}$ such that each
permutation equivalence class
of $\R^n$ is fully contained in a cone $C$ of $\mathcal{F}$. Such a fan
$\mathcal{F}$ defines a partition of $S_n$ because each maximal cone of the $S_n$-fan is contained in some cone $C \in \mathcal{F}$.

\begin{definition}
A {\em convex rank test} is a partition of the symmetric group
 $S_n$ which is induced by a coarsening of the $S_n$-fan.
We identify the fan with that rank test.
\end{definition}

We say that
two maximal cones, indexed by $\delta$ and $\delta'$, of the $S_n$-fan
{\em share a wall} if there exists an index $k$ such that
$\delta_k = \delta'_{k+1}$, $  \delta_{k+1} = \delta'_k$,
and $\delta_i = \delta'_i$ for $\,i \not\in \{k,k+1\}$. 
This condition means that
  the corresponding permutations $\delta$ and $\delta'$ differ
by an adjacent transposition.
To such an unordered pair $\{\delta,\delta'\}$,
we associate the following {\em (elementary) conditional independence (CI) statement}:
\begin{equation}
\label{CIStatement}
 \delta_k \perp \!\!\! \perp  \delta_{k+1} \,|\, \{\delta_1 , \ldots,  \delta_{k-1} \}.
 \end{equation}
The notation was coined by Dawid \cite{Dawid}, where it is used to formally describe conditional independence among sets of random variables; we will see the connection shortly.  
For $k=1$ we use the standard convention to abbreviate
$\, \delta_1 \perp \!\!\! \perp  \delta_{2} \,|\, \{\, \} \,$ by
$\, \delta_1 \perp \!\!\! \perp  \delta_{2}$.

\begin{example} \label{twenty-two}
For $n=3$ there are
$40$ pre-convex rank tests (Example \ref{bell}),
but only $22$ of them are convex rank tests.
The corresponding CI models are shown
in Figure 5.6 on page 108 in \cite{Studeny2005Probabilistic}.
\end{example}

The formula
(\ref{CIStatement})  defines a map from the set of walls of the $S_n$-fan
 onto the set 
$$
\mathcal{T}_n \,\, :=  \,\, \bigl\{
\, i \perp \!\!\! \perp j \,|\, K \,: \, K \subseteq [n] \backslash \{i,j\} \bigr\}. $$
of all elementary CI statements.
In this manner, each wall of the $S_n$-fan is labeled by a CI statement.
The map from walls to CI statements is not injective; 
there are $(n-k-1)!(k-1)!$ walls which are labeled by  \eqref{CIStatement}.  

The $S_n$-fan is the normal fan \cite{Ziegler1995}
of the {\em permutohedron} ${\bf P}_n$, which is the $(n-1)$-dimensional convex hull of
the vectors $(\rho_1,\ldots,\rho_n) \in \R^n$, where $\rho$ runs
over all rank vectors of permutations in $S_n$.
Each edge of ${\bf P}_n$ joins two permutations if they differ by an adjacent transposition.  In other words, each edge corresponds to a wall and is thus labeled by a CI statement.
A collection of parallel edges of ${\bf
P}_n$ that are perpendicular to a given hyperplane $\{x_i=x_j\}$ corresponds to the set of
CI statements $i \indep j |K$, where $K$ ranges over all subsets of
$[n] \backslash \{i,j\}$. 

The two-dimensional faces of ${\bf P}_n$
are squares and regular hexagons, and two edges of  ${\bf P}_n$ have
the same label in $\mathcal{T}_n$ if, but not only if, they are
opposite edges of a square.
Figure 2(c) depicts the subset of ${\bf P}_5$ in which the last two coordinates of $u \in \R^n$ are less than or equal to all other coordinates.  It consists of two copies of the hexagon in 2(a), with the final two entries
 of the descent vector either $4|5$ (in the top hexagon) or $5|4$ (in the bottom hexagon).
All vertical edges are labeled by the CI statement $4 \indep 5 | \{1,2,3\}$.

\begin{figure}[thb]\label{UpDown}
\[
\begin{array}{ccc}
 \begin{xy}<15mm,0cm>:
(-.9,1.3) *+!{++};
(-1.5,-.8) *+!{+-};
(1.5,.8) *+!{-+};
(.9,-1.4) *+!{--};
(-.5,\sinpioverthree)  ="123"  *+!DR{1|2|3} *{\bullet};
(.5 ,\sinpioverthree)  ="132"  *+!DL{1|3|2} *{\bullet};
(1,0)                  ="312"  *+!L{3|1|2}  *{\bullet};
(.5,-\sinpioverthree)  ="321"  *+!UL{3|2|1} *{\bullet};
(-.5,-\sinpioverthree) ="231"  *+!UR{2|3|1} *{\bullet};
(-1,0)                 ="213"  *+!R{2|1|3}  *{\bullet};
   "123";"132" **@{.};
   "132";"312" **@{-};
   "312";"321" **@{.};
   "321";"231" **@{.};
   "231";"213" **@{-};
   "213";"123" **@{.};
(\sinpioverthree, 0.4)  *+!{1 \indep 3 | \emptyset} ;
(-\sinpioverthree, -0.5) *+!{1 \indep 3 | \{ 2\} } ;
 \end{xy}
&
\begin{xy}
   <15mm,0mm>:
(0,0)  ="origin" ;
(\sinpioverthree, 0.5)  ="uprt"  *+!DL{1 \indep 3 | \emptyset} ;
(\sinpioverthree, -0.5)  ="dwnrt"  ;
(-\sinpioverthree, 0.5)  ="upl" ;
(-\sinpioverthree, -0.5)  ="dwnl" *+!UR{1 \indep 3 | \{ 2\} } ;
(0, 1)  ="up" ;
(0,-1)  ="dwn" ;
   "origin";"uprt" **@{--};
   "origin";"dwnrt" **@{-};
   "origin";"upl" **@{-};
   "origin";"dwnl" **@{--};
   "origin";"up" **@{-};
   "origin";"dwn" **@{-};
 \end{xy}
&
\begin{xy}<19mm,0cm>:
(-.5,.43)  ="123B"   *{\bullet};
(.5 ,.43)  ="132B"   *{\bullet};
(1,0)      ="312B"   *{\bullet};
(.5,-.43)  ="321B"   *{\bullet};
(-.5,-.43) ="231B"   *{\bullet};
(-1,0)     ="213B"   *{\bullet};
   "123B";"132B" **@{.};
   "132B";"312B" **@{.};
   "312B";"321B" **@{-};
   "321B";"231B" **@{-};
   "231B";"213B" **@{-};
   "213B";"123B" **@{.};
(-.5,.83)  ="123T"   *{\bullet};
(.5 ,.83)  ="132T"   *{\bullet};
(1,.4)    ="312T"    *{\bullet};
(.5,-.03)  ="321T"  *{\bullet};
(-.5,-.03) ="231T"  *{\bullet};
(-1,.4)     ="213T"  *{\bullet};
   "123T";"132T" **@{-};
   "132T";"312T" **@{-};
   "312T";"321T" **@{-};
   "321T";"231T" **@{-};
   "231T";"213T" **@{-};
   "213T";"123T" **@{-};
   "123T";"123B" **@{.};
   "132T";"132B" **@{.};
   "312T";"312B" **@{-};
   "321T";"321B" **@{-};
   "231T";"231B" **@{-};
   "213T";"213B" **@{-};
(-.075,0.58) *+!{4 \indep  5 |\{1,2,3\}};
(0,1) *+!{*|\!*\!|\!*\!|4|5};
(0,-.65) *+!{*|\!*\!|\!*\!|5|4};
 \end{xy}\\
\mbox{\bf (a)} & \mbox{\bf (b)}  & \mbox{\bf (c)}
\end{array}
\]
\caption{{\bf (a)} The permutohedron ${\bf P}_3$ and {\bf (b)} the $S_3$-fan projected to the plane.
The indicated rank test is up-down analysis. Each permutation is represented by
its descent vector $\delta = \delta_1 | \delta_2 | \delta_3$.  Missing walls of the $S_n$-fan, or solid edges of ${\bf P}_n$, are labelled by CI statements. 
 {\bf (c)} Edges of the permutohedron on opposite sides of a square (here, all vertical edges) are labelled by the same CI statement; hexagonal prisms such as the one pictured here appear in ${\bf P}_n$ for $n \geq 5$.}
\end{figure}
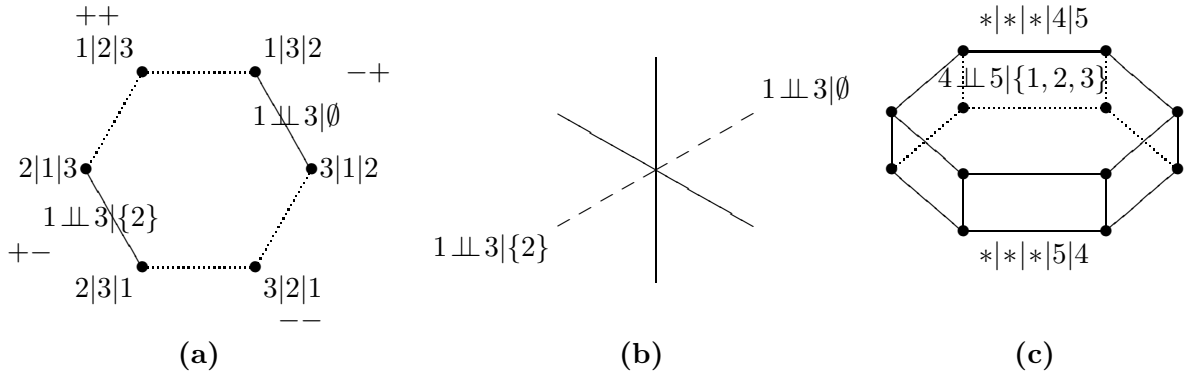

 Any convex rank test $\mathcal{F}$ is characterized
by the collection of walls $\{\delta,\delta'\}$ that are removed
when passing from the $S_n$-fan
to $\mathcal{F}$.
  So, from (\ref{CIStatement}), any convex rank test
$\mathcal{F}$ maps to a set  $\mathcal{M}_\mathcal{F}\,$
of CI statements corresponding to missing walls, or a set $\mathbf{M}_\mathcal{F}$ 
of edges of the permutohedron.
For example, if $\mathcal{F}$ is the fan obtained by removing the two dashed rays in 
Figure 2 (b) then the corresponding set
of CI statements is $\,\mathcal{M}_\mathcal{F} \,= \,
\bigl\{ 1 \indep 3 | \emptyset , \,1 \indep 3 | \{ 2\}  \bigr\}$.

Conditional independence statements \cite{Dawid} describe the dependence relationship among random variables.  
 A {\em semigraphoid}  is a set $\mathcal{M}$ of {\em general} conditional independence statements satisfying certain properties  \cite{Pearl}.
   These general conditional independence statements, in contrast to the elementary CI statements already introduced, can take subsets of $[n]$ in their first two arguments.  The conditions are, for $X,Y,Z$ pairwise disjoint subsets of $[n]$, 
\begin{eqnarray*}
{\rm (SG 1)} &  X \perp \!\!\! \perp Y \, |\, Z  \in \mathcal{M} \implies Y \perp \!\!\! \perp X \, |\, Z  \in \mathcal{M}\\
{\rm (SG 2)} &  X \perp \!\!\! \perp Y \, |\, Z  \in \mathcal{M} \mbox{ and } U \subset X \implies U \perp \!\!\! \perp Y \, |\, Z  \in \mathcal{M}\\
{\rm (SG 3)} &  X \perp \!\!\! \perp Y \, |\, Z  \in \mathcal{M} \mbox{ and } U \subset X \implies X \perp \!\!\! \perp Y \, |\, (U \cup Z)  \in \mathcal{M}\\
{\rm (SG 4)} &  X \perp \!\!\! \perp Y \, |\, Z  \in \mathcal{M} \mbox{ and } X \perp \!\!\! \perp W \,|\, (Y \cup Z) \implies X \perp \!\!\! \perp (W \cup Y) \, |\, Z  \in \mathcal{M}.
\end{eqnarray*}
It was shown by Studen\'{y} \cite{Studeny1990} that these are not a complete set of axioms for probabilistic conditional independence, although they are true of any probabilistic model.  A semigraphoid is determined by its {\em trace} among statements of the form $ i \perp \!\!\! \perp j \, |\, K $  where $i$ and $j$ are singletons.  Namely, $I \indep J | K$ holds if and only if $i \indep j |L$ for all $i \in I, j \in J$ and $L$ such that $K \subseteq L \subseteq (I \cup J \cup K) \setminus ij$; see \cite{Matus1992Equivalence}.
Casting the semigraphoid axiom in terms of the trace, we say that 
a subset $\mathcal{M} $ of $\mathcal{T}_n$ is a {\em semigraphoid} if 
$\, i \perp \!\!\! \perp j \, |\, K \in  \mathcal{M}\,$ implies
$\, j \perp \!\!\! \perp i \, |\, K \in  \mathcal{M}\,$ and
the following axiom holds:
\begin{eqnarray*}
{\rm (SG)} &&  \qquad \quad i \perp \!\!\! \perp j \, |\, K \cup {\ell}\, \in \mathcal{M}
\quad \, \mbox{and} \,\quad
i \perp \!\!\! \perp \ell \, |\, K\, \in \mathcal{M} \\
&& \!\!\! \mbox{implies } \,\,\,\,\,
i \perp \!\!\! \perp j \,|\, K \in \mathcal{M}
 \qquad \mbox{and } \quad i \perp \!\!\! \perp \ell \,| \, K \!\cup\! j  \in \mathcal{M}.
\end{eqnarray*}
This axiom is stated in \cite{Matus2004, Studeny2005Probabilistic}.
Our first result is that semigraphoids and convex rank tests are the same combinatorial object:

\begin{theorem} \label{fantheorem}
The map $\mathcal{F} \mapsto \mathcal{M}_\mathcal{F}$ is
a bijection between convex rank tests and  semigraphoids.
\end{theorem}

Before presenting the proof of this theorem, we shall discuss an example. 

\begin{example}[Up-down analysis] \label{ex.updwn} 
Let $\mathcal{F}$ denote the convex rank test called
 up-down analysis \cite{Willbrand2005}. In this test, each permutation
$\pi \in S_n$ is mapped to the sign vector of its first differences,
or, equivalently, its descent set. Thus this test is the natural map
$\,\tau : S_n \rightarrow \{-,+\}^{n-1}$.
The corresponding semigraphoid $ \mathcal{M}_\mathcal{F}$
consists of all CI statements $\, i \perp \!\!\! \perp j  \,| \, K \,$
where $\, | i -j | \geq 2 $.

This convex rank test is visualized in Figure 2(a,b) for $n=3$.  Permutations are in the same class (have the same sign pattern) if they are connected by a solid edge; there are four classes.  In the $S_3$-fan, the two missing walls are labeled by conditional independence statements as defined in (\ref
{CIStatement}). For $n=4$ the up-down analysis test  $\mathcal{F}$ 
is depicted in Figure 3.
The double edges correspond to the twelve CI statements
in $\mathcal{M}_\mathcal{F}$. There are eight classes; e.g.,
the class $\{3|4|1|2,3|1|4|2,1|3|4|2,1|3|2|4,3|1|2|4\}$ consists
of the five permutations in $S_4$ which have the up-down pattern  $(-,+,-)$.
\end{example}

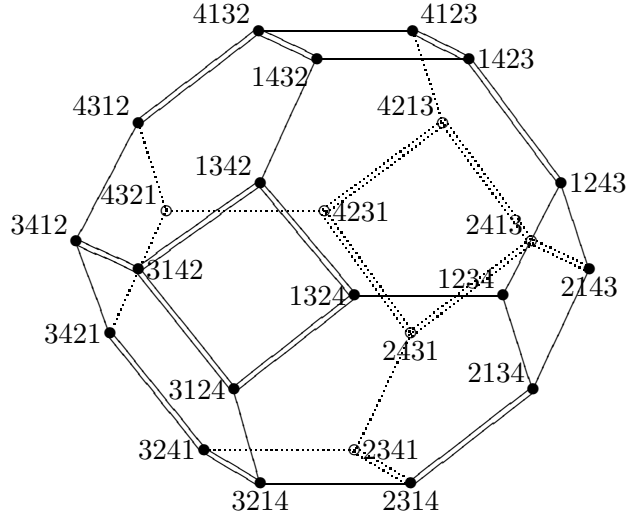
\begin{figure}[htb]\label{UpDown4}
\[
 \begin{xy}<25mm,0cm>:
(1,0)  ="3214"  *+!U{3214} *{\bullet};
(1.8,0) ="2314"  *+!U{2314} *{\bullet};
(.7,.18)  ="3241"  *+!R{3241} *{\bullet};
(1.5,.18)  ="2341"  *+!L{2341} *{\circ}; 
(.86,.5)  ="3124"  *+!R{3124} *{\bullet};
(2.45,.5)  ="2134"  *+!DR{2134} *{\bullet};
(.2,.8)  ="3421"  *+!R{3421} *{\bullet};
(1.8 ,.8)  ="2431"  *+!U{2431} *{\circ}; 
(1.5,1)  ="1324"  *+!R{1324} *{\bullet};
(2.29,1)  ="1234"  *+!DR{1234} *{\bullet};
(.35,1.14)  ="3142"  *+!L{3142} *{\bullet};
(2.75,1.14)  ="2143"  *+!U{2143} *{\bullet};
(.02,1.29)  ="3412"  *+!DR{3412} *{\bullet};
(2.44,1.29)  ="2413"  *+!DR{2413} *{\circ}; 
(.5,1.45)  ="4321"  *+!DR{4321} *{\circ}; 
(1.34,1.45)  ="4231"  *+!L{4231} *{\circ}; 
(1,1.6)  ="1342"  *+!DR{1342} *{\bullet};
(2.6,1.6)  ="1243"  *+!L{1243} *{\bullet};
(.35,1.92)  ="4312"  *+!DR{4312} *{\bullet};
(1.97,1.92)  ="4213"  *+!DR{4213} *{\circ}; 
(1.3,2.26)  ="1432"  *+!UR{1432} *{\bullet};
(2.11,2.26)  ="1423"  *+!L{1423} *{\bullet};
(.99,2.41)  ="4132"  *+!DR{4132} *{\bullet};
(1.81,2.41)  ="4123"  *+!DL{4123} *{\bullet};
"3214";"2314" **@{-}; 
"3241";"2341" **@{.}; 
"3241";"3214" **@{=}; 
"2341";"2314" **@{:}; 
"2134";"1234" **@{-}; 
"2143";"1243" **@{-}; 
"1234";"1243" **@{-}; 
"2134";"2143" **@{-}; 
"4132";"4123" **@{-}; 
"1432";"1423" **@{-}; 
"1432";"4132" **@{=}; 
"4123";"1423" **@{=}; 
"4312";"3412" **@{-}; 
"4321";"3421" **@{.}; 
"4312";"4321" **@{.}; 
"3412";"3421" **@{-}; 
"4213";"2413" **@{:}; 
"4231";"2431" **@{:}; 
"4213";"4231" **@{:}; 
"2413";"2431" **@{:}; 
"1342";"1324" **@{=}; 
"3142";"3124" **@{=}; 
"1342";"3142" **@{=}; 
"1324";"3124" **@{=}; 
"2314";"2134" **@{=}; 
"3124";"3214" **@{-}; 
"3421";"3241" **@{=}; 
"3412";"3142" **@{=}; 
"1324";"1234" **@{-}; 
"1432";"1342" **@{-}; 
"4312";"4132" **@{=}; 
"1423";"1243" **@{=}; 
"2341";"2431" **@{.}; 
"4321";"4231" **@{.}; 
"2413";"2143" **@{:}; 
"4123";"4213" **@{.}; 
\end{xy}
\]
\caption{The permutohedron ${\bf P}_4$
with vertices marked by descent vectors $\delta$ (bars $|$ omitted).
The convex rank test indicated by the double edges is up-down analysis.}
\end{figure}

Our proof of Theorem \ref{fantheorem} rests on translating
the semigraphoid axiom (SG) into
geometric statements about edges of the
permutohedron.
Recall that a semigraphoid $\mathcal{M}$ can be identified with the set $\mathbf{M}$ of
edges of the permutohedron whose CI statement labels are those of $\mathcal{M}$.

\begin{observation} \label{obs:SqHexAxioms}
A set $\mathbf{M}$ of edges of the permutohedron ${\bf P}_n$ is a
semigraphoid if and only if the set $\mathbf{M}$ satisfies the following
two geometric axioms: \\
{\bf Square axiom:} Whenever an edge of a square is in
$\mathbf{M}$, then the opposite edge is also
in $\mathbf{M}$. \\ 
\[
\begin{xy}<5mm,0cm>:
(-.7,.7)  ="TL"   *{\bullet};
(.7,.7)  ="TR"   *{\bullet};
(-.7,-.7)  ="BL"   *{\bullet};
(.7,-.7)  ="BR"   *{\bullet};
   "TL";"BL" **@{-};
   "TR";"BR" **@{.};
   "BR";"BL" **@{.};
   "TL";"TR" **@{.};
\end{xy} \quad
\implies \quad
\begin{xy}<5mm,0cm>:
(-.7,.7)  ="TL"   *{\bullet};
(.7,.7)  ="TR"   *{\bullet};
(-.7,-.7)  ="BL"   *{\bullet};
(.7,-.7)  ="BR"   *{\bullet};
   "TL";"BL" **@{-};
   "TR";"BR" **@{-};
   "BR";"BL" **@{.};
   "TL";"TR" **@{.};
\end{xy}
\]
{\bf Hexagon axiom:} Whenever two ad\-ja\-cent edg\-es of a hexagon
are in $\mathbf{M}$, then the two opposite edges of that hexagon
are also in $\mathbf{M}$.
\[
\begin{xy}<5mm,0cm>:
(-.5,\sinpioverthree)  ="123"   *{\bullet};
(.5 ,\sinpioverthree)  ="132"   *{\bullet};
(1,0)                  ="312"   *{\bullet};
(.5,-\sinpioverthree)  ="321"   *{\bullet};
(-.5,-\sinpioverthree) ="231"   *{\bullet};
(-1,0)                 ="213"   *{\bullet};
   "123";"132" **@{-};
   "132";"312" **@{.};
   "312";"321" **@{.};
   "321";"231" **@{.};
   "231";"213" **@{.};
   "213";"123" **@{-};
\end{xy} \quad
\implies \quad
\begin{xy}<5mm,0cm>:
(-.5,\sinpioverthree)  ="123"   *{\bullet};
(.5 ,\sinpioverthree)  ="132"   *{\bullet};
(1,0)                  ="312"   *{\bullet};
(.5,-\sinpioverthree)  ="321"   *{\bullet};
(-.5,-\sinpioverthree) ="231"   *{\bullet};
(-1,0)                 ="213"   *{\bullet};
   "123";"132" **@{-};
   "132";"312" **@{.};
   "312";"321" **@{-};
   "321";"231" **@{-};
   "231";"213" **@{.};
   "213";"123" **@{-};
\end{xy}
\]
\end{observation}

Let $\mathbf{M}$ be the subgraph of the edge graph of ${\bf P}_n$
defined by the statements in $\mathcal{M}$; that is, $\mathbf{M}$ consists of edges whose
labels are in $\mathcal{M}$.  Each class of the rank test defined by
$\mathcal{M}$ consists of the permutations in some connected
component of $\mathbf{M}$. We regard a path from
$\delta$ to $\delta'$ on ${\bf P}_n$ as a word $\sigma^{(1)} \cdots
\sigma^{(l)}$ in the free associative algebra $\mathcal{A}$
generated by the adjacent transpositions of $[n]$.  For example, the
transposition $\sigma_{23} := (23)$ gives the path from $\delta$ to
$\delta'=\sigma_{23} \delta = \delta_1 | \delta_3 | \delta_2 | \delta_4
| \dots  | \delta_n$.  The following relations in $\mathcal{A}$
define a presentation of the group algebra of $S_n$ as a quotient of $\mathcal{A}$:

\[\begin{array}{ccc}
{\rm (BS)} & \; \sigma_{i, i+1}\cdot \sigma_{i+k+1, i+k+2}\, -\, \sigma_{i+k+1, i+k+2} \cdot \sigma_{i, i+1}, & \\
{\rm (BH)} & \; \sigma_{i, i+1} \cdot \sigma_{i+1, i+2} \cdot \sigma_{i, i+1} - \sigma_{i+1, i+2}\cdot
 \sigma_{i, i+1} \sigma_{i+1, i+2}, & \qquad \mbox{   and} 
\\
{\rm (BN)} & \;  \sigma_{i, i+1}^2 -1, &
\end{array}\]
where suitable $i$ and $k$ vary over $[n]$.  
The first two are the \emph{braid relations}, and the third
represents the idempotency of each transposition.

Now, we regard these relations as properties of a
set of edges of ${\bf P}_n$, by identifying a word and a
permutation $\delta$ with the set of edges that comprise the
corresponding path in ${\bf P}_n$.  For example, a set satisfying
(BS) is one such that, starting from any $\delta$, the edges of the
path $\sigma_{i, i+1} \sigma_{i+k+1, i+k+2}$ are in the set if and
only if the edges of the path $ \sigma_{i+k+1, i+k+2} \sigma_{i,
i+1}$ are in the set. Note then, that (BS) is the square axiom,
and (BH) is a weaker version of the hexagon axiom of semigraphoids.  That is, implications in
either direction hold in a semigraphoid.  However, (BN) holds only
directionally in a semigraphoid: if an edge lies
in the semigraphoid, then its two vertices are in the same class;
but the empty path at some vertex $\delta$ certainly does not imply
the presence of all incident edges in the semigraphoid.  Thus, for
a semigraphoid, (BS) and (BH) hold, but (BN) must be replaced with
the directional version

\vspace{3mm}

\qquad \qquad \qquad  ${\rm (BN')} \;\; \
\qquad \sigma_{i, i+1}^2 \rightarrow 1.$
\vspace{3mm}

\noindent
We now consider a path $p$ from $\delta$ to $\delta'$ in a semigraphoid. 
Here is a crucial lemma for our proof:

\begin{lemma} \label{lem.allshortestpaths}
Suppose that $\mathcal{M}$ is a semigraphoid.
If $\delta$ and $\delta'$ lie in the same class of $\mathcal{M}$, then so do all shortest paths on ${\bf P}_n$ between them.
\end{lemma}

The lemma in turn depends on the following version of a classical result due to
Jacques Tits. This result, which can be found in~\cite[p.~49-51]{Brown1989}),
essentially states that  the relations
(BS),(BH),(BN) form a Gr\"obner basis for the
two-sided ideals they generate in $\mathcal{A}$.

\begin{theorem}[Tits \cite{Tits1968Problem}]
\label{thm:Tits}
Let $p$ and $q$ be words representing paths on $\mathbf{P}_n$.
\begin{itemize}
\item[(1)] A word $p$ is (BS),(BH),(BN)-reduced if and only if it is (BS),(BH),(BN')-reduced.
\item[(2)] If $p$ and $q$ are reduced, then they represent the same element of 
the symmetric group $S_n$ if and only if $p$ can be transformed to $q$ by the the application of (BS) and (BH) only.
\end{itemize}
\end{theorem}

\begin{proof}[Proof of Lemma \ref{lem.allshortestpaths}]
Theorem \ref{thm:Tits} (1) says that if there is any path connecting $\delta$ and $\delta'$, then there is a shortest path connecting them.  Thus if $\delta$ and $\delta'$ lie in the same class of $\mathcal{M}$, some shortest path $\delta \rightarrow \delta'$ also lies in that class.  Now (2) says that if $p$ and $q$ are both shortest paths, then $q$ can be obtained from $p$ by application of only the square and hexagon axioms, (BS) and (BH).  Thus if any shortest path $\delta \rightarrow \delta'$ lies in the class of $\mathcal{M}$ containing them both, so do all other shortest paths connecting them.
\end{proof}

We need one lemma to deal with intersections of nonmaximal cones.
Denote by $\isfaceof$ the transitive relation ``is a face of''
and write $F_w(C)$ for the face of a cone $C$
at which $w$ is minimized.

\begin{lemma} \label{lem:reduction}
If the intersection of two cones $C_1$ and $C_2$ is a face of
both, then the intersection of any faces $D \isfaceof C_1$ and $E \isfaceof C_2$ is a face of both.
\end{lemma}
\begin{proof}
By transitivity of $\isfaceof$ and the hypothesis it suffices to show $D \cap E \isfaceof C_1 \cap C_2$.  Since $D \isfaceof C_1$, there exists a 
linear functional $w$ such that the face $F_w(C_1)$ equals $D$ and $C_1 \cap C_2 \subset C_1 \subset H_w^+$.
Then $F_w(C_1 \cap C_2)=D \cap C_2$  so $D \cap C_2 \isfaceof C_1 \cap C_2$.
Similarly, $E \cap C_1 \isfaceof C_1 \cap C_2$.
Then since the intersection of any two faces of $C_1 \cap C_2$ is also
a face, $D \cap E \isfaceof C_1 \cap C_2$ as desired.
\end{proof}

\begin{proof}[Proof of Theorem \ref{fantheorem}]
Both semigraphoids and convex rank tests can be regarded as sets of edges of ${\bf P}_n$.
 We first
show that a semigraphoid satisfies (PC). Consider $\delta, \delta'$
in the same class $C$ of a semigraphoid, and let $\delta'' \in
\mathcal{L}(\delta \cap \delta')$.  Further, let $p$ be a shortest path from
$\delta$ to $\delta''$ (so, $p \delta = \delta''$), and let $q$ be a
shortest path from $\delta''$ to $\delta'$.  We claim that $qp$ is a
shortest path from $\delta$ to $\delta'$, and thus
$\delta'' \in C$ by Lemma
\ref{lem.allshortestpaths}.  Suppose $qp$ is not a
shortest path.  Then, we can obtain a shorter path in the
semigraphoid by some sequence of substitutions according to (BS),
(BH), and (BN').  Only (BN') decreases the length of a
path, so the sequence must involve (BN').  Therefore, there is some $i$, $j$ in $[n]$, such that their positions
relative to each other are reversed twice in $qp$. But $p$ and $q$
are shortest paths, hence  one reversal occurs in each of $p$ and
$q$. Then $\delta$ and $\delta'$ agree on whether $i>j$ or
$j>i$, but the reverse holds in $\delta''$, contradicting
$\delta'' \in \mathcal{L}(\delta \cap \delta')$.  Thus every
semigraphoid is a pre-convex rank test.

Now, we show that a semigraphoid corresponds to a fan.
We first argue that we may reduce to the case of two maximal cones, each coming from a class in the semigraphoid, whose intersection is codimension one in both.  By Lemma \ref{lem:reduction}, we can consider maximal cones only.  Suppose two maximal cones $C_1$, $C_k$ have intersection $C_1 \cap C_k$ which is not codimension one.  Then there exists a sequence of maximal cones $C_1, C_2, \dots, C_k$ such that $C_i \cap C_{i+1}$ is codimension one, $C_1 \cap C_k \subset C_i \cap C_{i+1}$ for all $i = 1, \dots k-1$, and in fact $C_1 \cap C_k = C_1 \cap C_2 \cap \cdots \cap C_k$.  We have that $(C_i \cap C_{i+1}) \cap (C_{i+1} \cap C_{i+2})$ is a face of $C_{i+1}$ and $C_{i+2}$ by Lemma \ref{lem:reduction}, and also is a face of $C_i$.  Thus $C_i \cap C_{i+1} \cap C_{i+2} \isfaceof C_i, C_{i+1}, C_{i+2}$; continuing in this manner, we eventually get that $C_1 \cap C_2 \cap \cdots \cap C_k \isfaceof C_1, C_k$ as required.

Consider the cone corresponding to a class $C$.  We need only
show that its codimension one intersection with another maximal cone is a shared face.
Since $C$ is a cone of a coarsening of the $S_n$-fan, each facet of $C$ lies in a hyperplane
$H =\{x_i=x_j\}$.
 Suppose a face of $C$ coincides with the hyperplane $H$
 and that $i>j$ in $C$.  A vertex $\delta$
borders $H$ if $i$ and $j$ are adjacent in $\delta$. We will
show that if $\delta,\delta' \in C$ border $H$, then their
reflections $\widehat{\delta} = \delta_1  | \dots |j|i| \dots
| \delta_n$ and $\widehat{\delta'}= \delta'_1 |  \dots |j|i| \dots |\delta'_n$
both lie in some class $C'$. Consider a `great circle' path between
$\delta$ and $\delta'$ which stays closest to $H$: all vertices
in the path have $i$ and $j$ separated by at most one position, and
no two consecutive vertices have $i$ and $j$ nonadjacent.  This is a
shortest path, so it lies in $C$, by Lemma
\ref{lem.allshortestpaths}. Using the
square and hexagon axioms (Observation \ref{obs:SqHexAxioms}), we
see that the reflection of the path across $H$ is a path
in the semigraphoid that connects $\widehat{\delta}$ to
$\widehat{\delta'}$ (Figure 3).  This shows that the intersection of $C$ and $C'$ is a face of both.  Thus a semigraphoid is a convex rank test.

Finally, if $\mathbf{M}$ is a set of edges of ${\bf P}_n$ representing
a convex rank test, then it is easy to show that $\mathbf{M}$
satisfies the square and hexagon axioms.
 \end{proof}
 
\begin{figure}[htb]\label{fig:reflection}
\[
 \begin{xy}<15mm,0cm>:
(0,0); 
p+ (\sinpioverthree, 0.5) *{\bullet}; 
p + (0,-1) *{\bullet} **@{.}; 
p + (-\sinpioverthree,-.5) *{\bullet} **@{-}; 
p + (-\sinpioverthree,+.5) *{\bullet} *+!UR{\widehat{\delta}} **@{-}; 
p + (0,1) *{\bullet} *+!DR{\delta} **@{.};
p+ (\sinpioverthree, 0.5) *{\bullet} **@{-};
p+ (\sinpioverthree, -0.5) *{\bullet} **@{-};
(\sinpioverthree, 0.5);
p + (1,0) **@{-};
(\sinpioverthree, -0.5);
p + (1,0) **@{-};
(2.73205081,0); 
p+ (\sinpioverthree, 0.5) *+!DL{\delta'} *{\bullet}; 
p + (0,-1) *{\bullet} *+!UL{\widehat{\delta'}} **@{.}; 
p + (-\sinpioverthree,-.5) *{\bullet} **@{-}; 
p + (-\sinpioverthree,+.5) *{\bullet} **@{-}; 
p + (0,1) *{\bullet} **@{.};
p+ (\sinpioverthree, 0.5) *{\bullet} **@{-};
p+ (\sinpioverthree, -0.5) *{\bullet} **@{-};
(-1.5,0);
p+(5.8,0) *+!UL{x_i = x_j}**@{--};
 \end{xy}\]
\caption{Reflecting a path across a hyperplane.}
\end{figure}
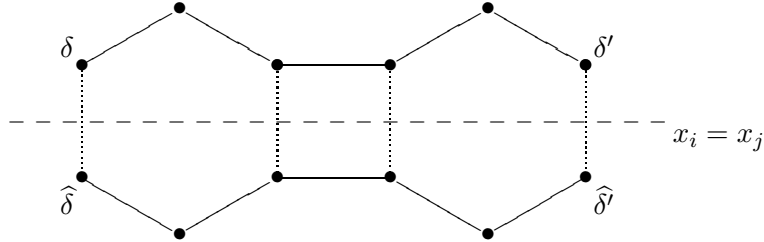

\section{The submodular cone}

In this section we focus on a subclass of the convex rank tests.
Let $2^{[n]}$ denote the collection of all
subsets of $[n] = \{1,2,\ldots,n\}$. Any real-valued function $\, w : 2^{[n]} \rightarrow \R \, $
defines a convex polytope $Q_w$ of dimension $\leq n-1$
as follows:
\begin{eqnarray*}  Q_w \,\,\, :=  &
\bigl\{ \, x \in \R^n \,: \,
x_1 + x_2 + \cdots + x_n = w([n]) \\
& \text{\  \,and } \sum\nolimits_{i \in I} x_i \leq w(I)\,\,
\hbox{for all} \,\, \emptyset\neq I \subseteq [n]   \,\bigr\}.
\end{eqnarray*}
A function $\, w : 2^{[n]} \rightarrow \R \, $ is called
 {\em submodular} if
$\,w(I) + w(J)\, \geq\, w(I \cap J) +  w(I \cup J)\,$
for $I,J \subseteq [n]$.  The {\em submodular cone}
is the cone ${\bf C}_n$ of all submodular functions $w :
2^{[n]} \rightarrow \R$.
  Working modulo its lineality space
   $\,{\bf C}_n \cap (-{\bf C}_n) $, we regard
   ${\bf C}_n$ as a pointed cone of dimension $2^n-n-1$.

Studying functions $w$ means that in considering the normal fan of a polytope $Q_w$, we want to retain information about non-binding inequalities that are just barely so, i.e.\ that hold with equality.  For this reason we define the {\em vector (normal) fan} \cite{BGS}.  
The indicator function of each $I \in 2^{[n]}$ defines a vector $e_I$ in the $1$-skeleton of the $S_n$-fan, understood modulo $e_{[n]}$; for example, these vectors for $n=3$ are $e_{001}, e_{010}, e_{100}, e_{011}, \dots, e_{111}$.  
A {\em vector fan} $\frkF$ is a collection of subsets of  $\{e_I: I \in 2^{[n]}\}$ such that $U, V \in \frkF$ implies $U \cap V \in \frkF$.
A vector fan defines a usual fan by taking the maximal cones of the fan to be the cones generated by the vector sets in the vector fan.  We say that a vector fan is \emph{complete} if its fan is.  A vector fan $\frkF$ \emph{coarsens} another vector fan $\frkG$ if for all $U \in \frkG$, there exists $V \in \frkF$ with $U \subset V$. 

Given a function $w:2^{[n]} \rightarrow \R$, each $I\in 2^{[n]}$ defines an inequality $\sum_{i \in I} x_i \leq w_I$ appearing in the definition of $Q_w$; the vector normal fan tells us which of these inequalities holds with equality on some face of $Q_w$.  We define the {\em vector normal fan} of a function $w:2^{[n]} \rightarrow \R$ as the set $\{ \{e_I: I \in 2^{[n]}, \sum_{i \in I}x_i =w_I$ for all $x \in F\}$ for each face $F \in Q_w \}$. The vector normal fan of $w$ defines a fan which is the normal fan of $Q_w$ and retains additional information.

\begin{proposition} \label{prop:submodularnormal} 
A function $\,w: 2^{[n]} \rightarrow \R \, $
is submodular if and only if
the vector normal fan of $w$ is a coarsening of the vector $S_n$-fan.
\end{proposition}

\begin{example}
Let $w_1 = w_2 =w_3 =1, w_{12}=w_{13}=w_{23}=w_{123}=3$. 
The polytope $Q_w$ is the point $(1,1,1)$ but the function $w$ is not submodular. The vector normal fan $\frkF$ of $w$ is $\{\{e_{001},e_{010},e_{100}\}\}$ and the normal fan is all of $\R^3 / (1,1,1)$.  $\frkF$ does not coarsen the $S_n$-fan since, for example, $e_{110}$ is not contained in any set in $\frkF$.

However, if we change $w$ slightly to define the same $Q_w$ but with the inequalities corresponding to $011,101$, and $110$ also holding with equality, e.g.
$w_1 = w_2 =w_3 =1, w_{12}=w_{13}=w_{23}=2$, and $w_{123}=3$, the resulting vector normal fan of $w$ is a coarsening of the (vector) $S_n$-fan.
\end{example}

\begin{proof}
We show only the if direction
of Proposition \ref{prop:submodularnormal}.
  Suppose $w$ is not submodular.  Then there exist $I,J \subset 2^{[n]}$ such that\[
w_I + w_J < w_{I \cap J} + w_{I \cup J}
\] 
We also have that
\begin{eqnarray*}
\sum_{i \in I \cup J} x_i + \sum_{i \in I \cap J} x_i & = & \sum_{i \in I} x_i + \sum_{i \in J} x_i\\
&\leq & w_I + w_J < w_{I \cap J} + w_{I \cup J}
\end{eqnarray*}
So $\sum_{i \in I \cup J} x_i < w_{I \cup J} + (w_{I \cap J}- \sum_{i \in I \cap J} x_i)$ and similarly $\sum_{i \in I \cap J} x_i < w_{I \cap J} + (w_{I \cup J}- \sum_{i \in I \cup J} x_i)$, so that at most one of the inequalities corresponding to $I \cup J$ and $I \cap J$ can hold with equality at any point of $Q_w$.  Then any set in the vector normal fan of $w$ either fails to
contain $e_{I \cap J}$ or fails to contain $e_{I \cup J}$.
\end{proof}

Proposition \ref{prop:submodularnormal} can be paraphrased as follows:
 the function $w$ is submodular if and only if
the optimal solution of
$$
\mbox{maximize $u \cdot x$ subject to $x \in Q_w$}
$$
depends only on the permutation equivalence class
of $u$.
Thus, solving this linear programming problem
constitutes a convex rank test.  Any such test is called a
{\em submodular rank test}.

A convex polytope is a {\em (Minkowski) summand}
of another polytope if the normal fan of the latter
refines the normal fan of the former. The
polytope $Q_w$ that represents a submodular rank test
is a  summand of the permutohedron
${\bf P}_n$.

\begin{theorem}
The following combinatorial objects are equivalent for any positive integer~$n$: \\
\noindent $1.$ submodular rank tests, \hfill \break
\noindent $2.$ summands of the permutohedron $\mathbf{P}_n$, \hfil \break
\noindent $3.$ structural conditional~independence~models \cite{Studeny2005Probabilistic}, \hfil \break
\noindent $4.$ faces of the submodular cone ${\bf C}_n$ in $\R^{2^n}$. 
\end{theorem}

\begin{proof}
We have 1$\iff$2 from Proposition \ref{prop:submodularnormal}, and
1$\iff$3 follows from \cite{Studeny2005Probabilistic}. Further, 1$\iff$4 
is a direct consequence of our definition of submodular rank tests.
\end{proof}

\begin{remark}
All $22$ convex rank tests for $n=3$ are submodular.
The submodular cone ${\bf C}_3$ is a
$4$-dimensional cone whose base is a
bipyramid. Its f-vector is $(1,5,9,6,1)$. The polytopes $Q_w$, as
$w$ ranges over representatives of the faces of ${\bf C}_3$,
are all the Minkowski summands of~${\bf P}_3$.
\end{remark}

\begin{proposition} \label{notsubmodular}
For $n \geq 4$, there exist convex rank tests that are not submodular rank tests.
Equivalently, there are fans that coarsen the $S_n$-fan
but are not the normal fan of any polytope.
\end{proposition}

\begin{proof}
This result is well-known. It is stated in Section 2.2.4 of \cite{Studeny2005Probabilistic} in the following form:  ``There exist semigraphoids that are not structural.''
\end{proof}

An interesting example which also proves Proposition \ref{notsubmodular}
is the following semigraphoid:
$$
\mathcal{M} \quad = \quad
\bigl\{
 2 \perp \!\!\! \perp 3 | \{1,4\},\,
 1 \perp \!\!\! \perp 4 |  \{2,3\}, \,  
 1 \perp \!\!\! \perp 2 | \emptyset,\,
 3 \perp \!\!\! \perp 4 |\emptyset \,\bigr\}.
 $$
The corresponding fan consists of unimodular cones, or, equivalently,
the posets $P_i$ representing this non-submodular convex rank test are all trees.
This example answers a question posed in the first version of \cite{PRW}.
A systematic method for showing that a semigraphoid is not submodular can be
found in \cite{counterexamples}. Results in that paper include
an example of a coarsest semigraphoid which is not submodular and a proof that the semigraphoid semigroup is not normal.

\begin{remark} \label{rmk1}
For $n=4$ there are $22108$ submodular rank tests, one for each face of the
$11$-dimensional cone ${\bf C}_4$.
The base of this submodular cone is a
$10$-dimensional polytope with
$f$-vector $
(1,37, 356,  $ $ 1596, 3985, 5980, 5560, 3212, 1128, 228, 24,1)$.  The $37 $
vertices of this polytope correspond to the maximal semigraphoids.
These come in seven
symmetry classes up to the $*$ involution (\ref{starinvolution}) and the $S_4$-action.
The types of maximal semigraphoids
for $n=4$ are displayed in the following table:

\bigskip
{\small
\noindent \begin{tabular}{lcccc}
Symmetry & \!\!\!\! No. \!\!\!\!\!\! & $i \indep j $  & $i \indep j | k$ & $i \indep j | \{k,l\}$ \\
\hline
$1 \times$ and $*$ & 2 & all & all & none \\
$4 \times$ and $*$ & 8 & all & \!\! \! all but $2 \indep 3 | 1, 1 \indep 3 | 2, 1 \indep 2 |3$ & $3 \indep 4|12, 2 \indep 4 |13, 1 \indep 4|23$\\
$6 \times$ incl. $*$ & 6 & all but $1\indep 2$& all but $1 \indep 2 | 3, 1 \indep 2 | 4$ & all but $1 \indep 2|34$\\
$4 \times$ and $*$ & 8 & all & $2 \indep 3 | 4, 2 \indep 4 | 3, 3 \indep 4 |2$ & $3 \indep 4|12, 2 \indep 4 |13, 2 \indep 3|14$\\
$1 \times$, self-$*$ & 1 & all & none & all \\
$6 \times$ incl. $*$ & 6 & all but $1\indep 2$& $2 \indep 3|1, 2\indep 4|1, 1 \indep 3|2, 1 \indep 4|2$ & all but $3 \indep 4|12$\\
$6 \times$ incl. $*$ & 6 & $3\indep 4$& all but $2 \indep 3 | 4, 2\indep 4 | 3, 1 \indep 4 | 3, 1 \indep 3|4$ & $1 \indep 2|34$\\
\hline
\end{tabular}
}
\end{remark}
\medskip

\begin{remark} \label{rmk2}
For $n=5$ there are
$117978$ coarsest submodular rank tests,
in $1319$ $S_5$ symmetry classes.
We confirmed this result of \cite{Studeny2000} with {\tt POLYMAKE} \cite{Gawrilow2000}.
\end{remark}

\medskip

We now define a class of submodular rank tests,
which we call {\em Minkowski sum of
simplices (MSS) tests}. 
Note that each subset $K$ of $[n]$
defines a submodular function $w_K$
by setting $w_K (I) = 1$ if $K \cap I $ is non-empty
and  $w_K(I) = 0$ if $K \cap I $ is empty.
The corresponding polytope
$Q_{w_K}$ is the simplex
$\Delta_K = {\rm conv} \{ e_k :  k \in K \}$.

Now consider an arbitrary subset  $\,\mathcal{K} = \{K_1,K_2,\ldots,K_r \}\,$
of $2^{[n]}$. It defines the submodular function
$\,w_{\mathcal{K}} =  w_{K_1} + w_{K_2} + \cdots + w_{K_r}$.
 The corresponding polytope is the Minkowski sum
$$ \Delta_\mathcal{K} \quad = \quad \Delta_{K_1} + \Delta_{K_2} + \cdots + \Delta_{K_r}. $$
 The associated MSS test $\tau_\mathcal{K}$ is defined as follows.
 Given $\rho \in S_n$, we compute the number of indices
 $j \in [r]$ such that $\,{\rm max}\{ \rho_k \,: \, k \in K_j \}\, = \,
\rho_i $,
 for each  $i \in [n]$.
 The signature $\tau_\mathcal{K}(\rho)$ is
the vector in $\N^n$ whose $i$th coordinate is that number. 
Few submodular rank tests are MSS tests:

\begin{remark} \label{rmk3}
For $n = 3$, there are $22$ submodular rank tests,
but only $15$ of them are MSS tests.
For $n=4$, there are $22108$ submodular rank tests,
but only  $1218$ of them are MSS tests.
\end{remark}

In light of Theorem \ref{fantheorem}, it is natural to ask
which semigraphoids correspond to an MSS test.
 Geometrically, we wish to know which edges of
 the permutohedron ${\bf P}_n$ are contracted
 when passing to the  polytope $Q_{w_{\mathcal{K}}}$.
 To be precise,  let  $\mathcal{M}_\mathcal{K}$ denote
 the semigraphoid derived from $\mathcal{F}_{w_{\mathcal{K}}}$ using the bijection in
Theorem \ref{fantheorem}. We then have the following result:

\begin{proposition} \label{CIsetfam}
The semigraphoid $\mathcal{M}_\mathcal{K}$
is the set of CI statements of the form
$\, i \perp \!\!\! \perp j \, |\, K \,$ where all
sets containing $\{i,j\}$
and contained in $\{i,j\} \cup [n] \backslash K \,$
are not in $\mathcal{K}$.
\end{proposition}

\begin{proof}
Consider two permutations  $\delta$ and $\delta'$
which are adjacent on the permutohedron ${\bf P}_n$,
and let $\,i \perp \!\!\! \perp j \, |\, K \,$ be the label of the edge
that connects $\delta$ and $\delta'$.
That CI statement is in $\mathcal{M}_\mathcal{K}$
if and only if $\delta$ and $\delta'$ are mapped
to the same vertex in $\Delta_{\mathcal{K}}$ if and only if
$\delta$ and $\delta'$ are mapped
to the same vertex in each simplex $\Delta_{K_l}$ 
for $l=1,2,\ldots,r$. For each $l$, this means that
the leftmost entry of the descent vector $\delta$ that lies in $K_l$
agrees with the leftmost entry of the other descent vector $\delta'$ that lies in $K_l$.
This condition is equivalent to 
$$ K_l \,\, \cap \,\, (\,K \,\cup \,\{i,j\} \,) \quad \not= \quad \{i,j\} \qquad \qquad
\hbox{for}\, \,\,\, l =1,2,\ldots,r .$$
Thus  $\,i \perp \!\!\! \perp j \, |\, K \,$ is in the semigraphoid
$\,\mathcal{M}_{\mathcal{K}}\,$ associated with the set family $\mathcal{K}$
 if and only if $\mathcal{K}$
contains no set whose intersection with
$\,K\,\cup \,\{i,j\} \,$ equals $\,\{i,j\}$.
This is precisely our claim.
\end{proof}

There is a natural involution $*$ on the set of all CI statements
which is defined as follows:
\begin{equation} \label{starinvolution}
 ( i \perp \!\!\! \perp j \,|\,  C)^* \quad := \quad
 i \perp \!\!\! \perp j \,|\,  [n]\backslash (C \cup \{i,j\}) .
\end{equation}
If $\mathcal{M}$ is any semigraphoid,
then the semigraphoid $\mathcal{M}^*$ is obtained by applying the involution $*$
to all the CI statements in the model $\mathcal{M}$.  This involution is referred to as {\em duality}
 in \cite{Matus1992Ascending}. In the {\em boolean lattice},
whose elements are the subsets of $[n]$, the involution corresponds to
switching the role of set intersection and set union.

The MSS test  $\tau_{\mathcal{K}}$ was defined above in terms of
 weight functions $w$. What follows is a similar construction for the
 duals of MSS tests.
Let $z_{\mathcal{K}}(J)=1$ for $J \in
\mathcal{K}$ and $z_{\mathcal{K}}(J) =0$ otherwise.  Then the function
$\,w^*: 2^{[n]} \rightarrow \R \, $ defined by
 $\,w_{\mathcal{K}}^* (I) := \sum_{J \subset I} z_{\mathcal{K}}(J)\,$ is supermodular. We set
\begin{eqnarray*}  Q_w^* & :=  \,
\bigl\{ \, x \in \R^n \,: \,
x_1 + x_2 + \cdots + x_n = w([n]) \\
& \text{\ \,and } \sum\nolimits_{i \in I} x_i \geq w(I)\,\,
\hbox{for all} \,\, \emptyset\neq I \subseteq [n]   \,\bigr\}.
\end{eqnarray*}
Then the equality $\,Q_{w_{\mathcal{K}}}^* \,=\, \Delta_K\,$ holds
for $\,\Delta_{\mathcal{K}} = \Delta_{\mathcal{K}_1} + \Delta_{K_2} + \cdots + \Delta_{K_r}$.
This equality is precisely the statement in Proposition 6.3 of
Postnikov's paper \cite{Postnikov2005}. 

\section{Graphical tests} \label{sec:graphical}

We have seen that semigraphoids are equivalent to convex rank tests.
We now explore the connection to graphical models.
 Let $G$ be a graph with vertex set $[n]$ and
$\mathcal{K}(G)$ the collection of all subsets $K \subseteq
[n]$ such that the induced subgraph of $G|_K$ is connected.
The {\em undirected graphical model}
(or {\em Markov random field}) derived from the graph $G$
is the set $\mathcal{M}^G$  of CI statements:
\begin{equation}
\label{noPath}
  \mathcal{M}^G \,\,\, = \,\,\,
\bigl\{\, i \perp \!\!\! \perp j \,|\,  C \,\, :\,\,
\mbox{the restriction of $\,G\,$ to}  \,\,\,  [n] \backslash C \,\,
\mbox{ contains no path from $i$ to $j$} \bigr\}.
\end{equation}

\begin{theorem} \label{Jasonslemma}
The set $\mathcal{M}^G$  of CI statements in the graphical model $G$ is equal to the
semigraphoid $\,\mathcal{M}_{\mathcal{K}(G)}\,$
associated with the family $\mathcal{K}(G)$
of connected induced subgraphs of $G$.
\end{theorem}

\begin{proof}
 The defining condition in (\ref{noPath}) is equivalent
 to saying that the restriction of $G$ to any node
 set containing $\{i,j\}$ and contained in $\,\{i,j\} \cup ([n] \backslash C)\,$
 is disconnected. With this observation, Theorem 
 \ref{Jasonslemma} follows directly from
 Proposition~\ref{CIsetfam}.
\end{proof}

The polytope $\Delta_G = \Delta_{\mathcal{K}(G)}$
associated with the graph $G$ 
is the {\em graph associahedron}. This is
a well-studied object in combinatorics
\cite{Postnikov2005,Carr2004}.
Carr and Devadoss \cite{Carr2004} showed that
$\Delta_G$ is a simple polytope whose
faces are in bijection with the tubings
of the graph $G$. Tubings are defined as follows.
Two subsets $A,B$ $\subset [n]$ are
\emph{compatible} for $G$
if one of the following conditions holds: $A\subset B$, $B\subset A$, or $A\cap B = \emptyset$, and there is no edge between any node in
$A$ and $B$. A {\em tubing} of the graph $G$
is a subset ${\bf T}$ of $2^{[n]}$ such that
any two elements of ${\bf T}$ are compatible.
The set of all tubings on $G$ is a simplicial complex;
it is dual to the face lattice of the simple polytope $\Delta_G$.

For any graph $G$ on $[n]$ we now have two convex rank tests.
First, there is the {\em graphical model rank test} $\,\tau_{\mathcal{K}(G)}$, 
which is the MSS test of the set family $\mathcal{K}(G)$. Second,
we have the {\em graphical tubing rank test} $\,\tau^*_{\mathcal{K}(G)}$, which
is the convex rank test associated with the semigraphoid $\,(\mathcal{M}^G)^*\,$
dual to $\,\mathcal{M}^G$. Explicitly, that dual semigraphoid is given by
\begin{equation}
\label{noPath2}
\!  (\mathcal{M}^G)^* \,\, = \,\,
\bigl\{\, i \perp \!\!\! \perp j \,|\,  C \,\, :\,
\mbox{the restriction of $\,G\,$ to}  \,\,\,   C  \cup \{i,j\} \,
\mbox{ contains no path from $i$ to $j$} \bigr\}.
\end{equation}

\begin{figure}[htb]\label{GM}
\[
 \begin{xy}<25mm,0cm>:
(1,0)  ="3214"  *+!U{3214} *{\bullet};
(1.8,0) ="2314"  *+!U{2314} *{\bullet};
(.7,.18)  ="3241"  *+!R{3241} *{\bullet};
(1.5,.18)  ="2341"  *+!L{2341} *{\circ}; 
(.86,.5)  ="3124"  *+!R{3124} *{\bullet};
(2.45,.5)  ="2134"  *+!DR{2134} *{\bullet};
(.2,.8)  ="3421"  *+!R{3421} *{\bullet};
(1.8 ,.8)  ="2431"  *+!U{2431} *{\circ}; 
(1.5,1)  ="1324"  *+!R{1324} *{\bullet};
(2.29,1)  ="1234"  *+!DR{1234} *{\bullet};
(.35,1.14)  ="3142"  *+!L{3142} *{\bullet};
(2.75,1.14)  ="2143"  *+!U{2143} *{\bullet};
(.02,1.29)  ="3412"  *+!DR{3412} *{\bullet};
(2.44,1.29)  ="2413"  *+!DR{2413} *{\circ}; 
(.5,1.45)  ="4321"  *+!DR{4321} *{\circ}; 
(1.34,1.45)  ="4231"  *+!L{4231} *{\circ}; 
(1,1.6)  ="1342"  *+!DR{1342} *{\bullet};
(2.6,1.6)  ="1243"  *+!L{1243} *{\bullet};
(.35,1.92)  ="4312"  *+!DR{4312} *{\bullet};
(1.97,1.92)  ="4213"  *+!DR{4213} *{\circ}; 
(1.3,2.26)  ="1432"  *+!UR{1432} *{\bullet};
(2.11,2.26)  ="1423"  *+!L{1423} *{\bullet};
(.99,2.41)  ="4132"  *+!DR{4132} *{\bullet};
(1.81,2.41)  ="4123"  *+!DL{4123} *{\bullet};
"3214";"2314" **@{-}; 
"3241";"2341" **@{.}; 
"3241";"3214" **@{=}; 
"2341";"2314" **@{:}; 
"2134";"1234" **@{-}; 
"2143";"1243" **@{-}; 
"1234";"1243" **@{-}; 
"2134";"2143" **@{-}; 
"4132";"4123" **@{-}; 
"1432";"1423" **@{-}; 
"1432";"4132" **@{*}; 
"4123";"1423" **@{*}; 
"4312";"3412" **@{-}; 
"4321";"3421" **@{.}; 
"4312";"4321" **@{.}; 
"3412";"3421" **@{-}; 
"4213";"2413" **@{o}; 
"4231";"2431" **@{o}; 
"4213";"4231" **@{:}; 
"2413";"2431" **@{:}; 
"1342";"1324" **@{=}; 
"3142";"3124" **@{=}; 
"1342";"3142" **@{*}; 
"1324";"3124" **@{*}; 
"2314";"2134" **@{=}; 
"3124";"3214" **@{-}; 
"3421";"3241" **@{=}; 
"3412";"3142" **@{=}; 
"3412";"3142" **@{*}; 
"1324";"1234" **@{-}; 
"1432";"1342" **@{-}; 
"4312";"4132" **@{*}; 
"1423";"1243" **@{*}; 
"2341";"2431" **@{.}; 
"4321";"4231" **@{.}; 
"2413";"2143" **@{:}; 
"2413";"2143" **@{o}; 
"4123";"4213" **@{.}; 
\end{xy}
\]
\caption{The permutohedron ${\bf P}_4$. Double edges indicate the MSS test $\tau_{\mathcal{K}(G)}$ where
$G$ is the $4$-chain. Edges with large dots
indicate the dual tubing test $\,\tau^*_{\mathcal{K}(G)}$.}
\end{figure}
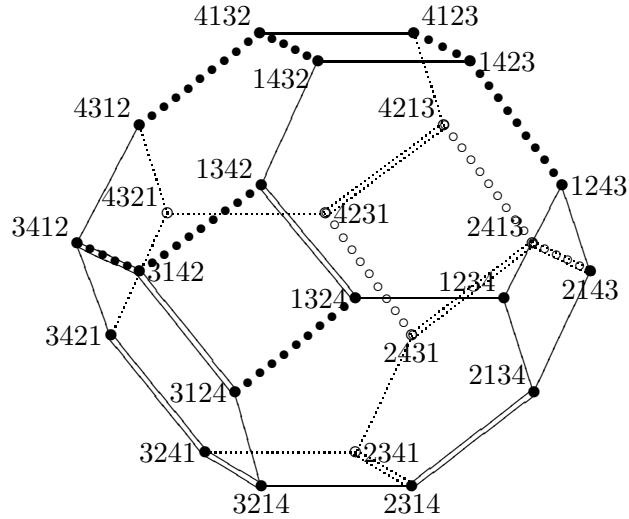

We summarize our discussion in the following theorem:

\begin{theorem} \label{maingraphical}
The following four
combinatorial objects are isomorphic for any graph $G$ on $[n]$: \hfill \break
\noindent $\bullet$ the graphical model rank test $\tau_{\mathcal{K}(G)}$, \hfill \break
\noindent $\bullet$ the graphical tubing rank test
$\tau^*_{\mathcal{K}(G)}$, \hfill \break
\noindent $\bullet$ the fan of the graph associahedron~$\Delta_G$, \hfill \break
\noindent $\bullet$ the simplicial complex of all tubings on $G$.
\end{theorem}

We note that when the graph $G$ is a path of length $n$, $\Delta_G$ is the
{\em associahedron}, and
when it is an $n$-cycle, $\Delta_G$ is the {\em cyclohedron}. The number of classes in 
either the MSS test $\tau_{\mathcal{K}(G)}$ or the
tubing test $\tau^*_{\mathcal{K}(G)}$ is the
{\em $G$-Catalan number} of \cite{Postnikov2005}.  This number is
the classical Catalan number
 $\frac{1}{n+1} {2n \choose n}$ for the associahedron test.
It equals ${2n-2 \choose n-1}$ for the cyclohedron test \cite{cyclohedron}.

\begin{example}
Let $n=4$ and let $G $ be
the $4$-chain $\, 1$---$2$---$3$---$4$. Then
$$
\begin{matrix}
  \mathcal{M}^G   \!\! & = & \bigl\{
  1 \perp \!\!\! \perp 3 \,|\,  24, & 
    1 \perp \!\!\! \perp 4 \,|\,  23, & 
    2 \perp \!\!\! \perp 4 \,|\,  13, & 
  1 \perp \!\!\! \perp 3 \,|\,  2, & 
1 \perp \!\!\! \perp 4 \,|\,  2, & 
1 \perp \!\!\! \perp 4\,|\,  3, & 
2 \perp \!\!\! \perp 4 \,|\,  3 \bigr\},\\ 
(\mathcal{M}^G)^* \!\! &  = & \bigl\{
1 \perp \!\!\! \perp 3 \, ,& 
    1 \perp \!\!\! \perp 4 \, , & 
    2 \perp \!\!\! \perp 4 \,,  & 
  1 \perp \!\!\! \perp 3 \,|\,  4, & 
1 \perp \!\!\! \perp 4 \,|\,  3, & 
1 \perp \!\!\! \perp 4\,|\,  2, & 
2 \perp \!\!\! \perp 4 \,|\,  1 \bigr\}.
\end{matrix}
$$
The corresponding tests $\tau_{\mathcal{K}(G)}$ and
$\tau_{\mathcal{K}(G)}^*$ are depicted in Figure \ref{GM}.
Note that contracting either class of marked edges on the
permutohedron in Figure \ref{GM} leads
to the $3$-dimensional associahedron $\Delta_G$.
The associahedron $\Delta_G$ is the Minkowski sum of 
the simplices $\Delta_K$ where $K$ runs over 
$$\mathcal{K}(G) \quad = \quad \bigl\{ \{1\}, 
\{2\}, \{3\}, \{4\}, \{1,2\}, \{2,3\}, \{3,4\}, \{1,2,3\}, \{2,3,4\}, \{1,2,3,4\} \bigr\} . $$
The $3$-dimensional simple polytope $\Delta_4$ has $14$ vertices,
one for each of the $14$ tubings of $G$. \qed
\end{example}

In our application of graphical rank tests,
we found it more natural to work with the tubing test
$\tau^*_{\mathcal{K}(G)} $ instead of the MSS test
$\tau_{\mathcal{K}(G)}$. We refer to our
 companion paper \cite{cyclohedron} which gives
a detailed discussion of the cyclohedron test and
its applications.
By the cyclohedron test we mean the tubing test  $\,\tau^*_{\mathcal{K}(G)}\, $  where
the graph $G$ is a cycle of length $n$.

\begin{figure}[htb]
\includegraphics[scale=0.5]{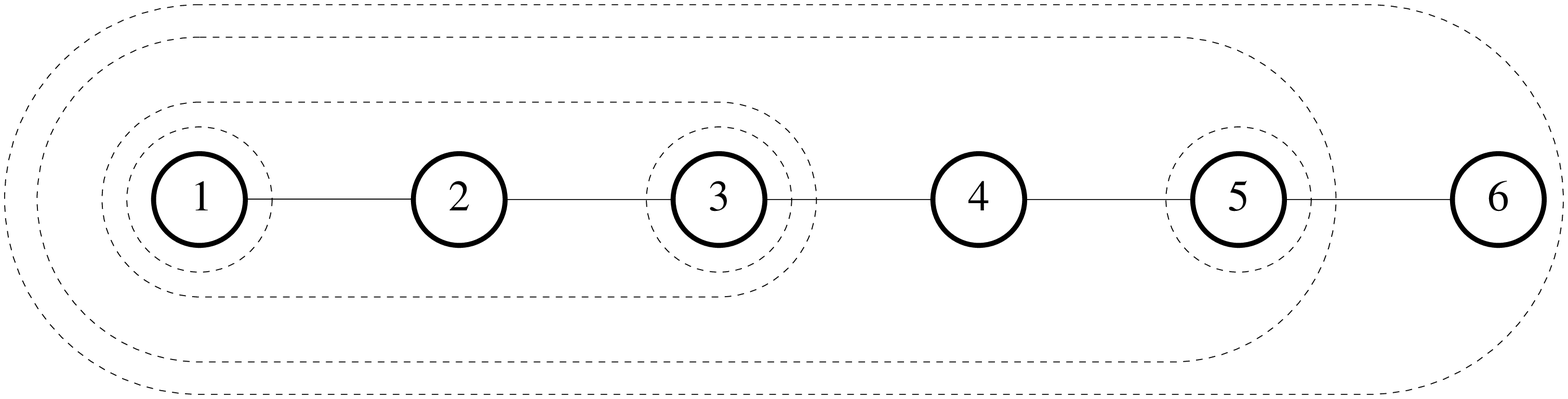}
\caption{Tubing of the $6$-chain.  Encircled regions indicate the sets $U_j.$}
\end{figure}

Applying the tubing test to a data vector $u \in \R^n$
 can be viewed as an iterative procedure for drawing a topographic map on the graph $G$.
Namely, we encircle the vertices of $G$ by sets $U_1, \dots, U_n$  
in the order $\delta_1, \delta_2, \dots, \delta_{n-1}$, with the following provision: if $\delta_i$ is next to be encircled and shares an edge with some vertex $j$ which has already been encircled by some $U_j$, then $U_i$ must also contain the circle $U_j$.
The result is  a collection $U$ of $n-1$ encircled sets $U_1,U_2,\ldots,U_{n-1}$,
and this unordered collection of sets 
is the signature of $v$. The height $h_i$ of the $i$-th node in the
topographic map for $v$ is the number of sets $U_j$ which contain $i$. 
We can identify the signature $U$ with the
{\em height vector} $h = (h_1,h_2,\ldots,h_n)$, since
$U$ can be recovered uniquely from the vector $h$.
The map $u \mapsto h(u)$ can be interpreted
as a {\em smoothing of the data}.  Figure 6 displays the topographic map when the data vector is 
$\,u=(2.1,0.3 ,1.8,,2.0,1.1,0.1)$.
Here $G$ is the 
$6$-chain $\, 1$---$2$---$3$---$4$---$5$---$6$.
and the descent vector of $u$  equals $\delta=(1|5|3|2|4|6)$.

\section{On counting linear extensions}

In this paper, we have  introduced a hierarchy of rank tests, which range from
pre-convex to graphical. Convex rank tests are applied to data
vectors $u \in \R^n$, or permutations $\pi \in S_n$, and determine
their cones in a fan $\mathcal{F}$ which coarsens the $S_n$-fan.
The significance of a data vector in such a test is measured
by a certain p-value, whose precise derivation is described in 
\cite{cyclohedron}. Computation of that p-value
rests on our ability to compute the quantity $\,|\, \tau^{-1} \bigl(
\tau(\pi) \bigr)\,|$, which is the number of permutations in the maximal cone 
of $\mathcal{F}$ corresponding to $\pi$. 
Recall that the cones  of a
convex rank test are indexed by posets $P_1,P_2,\ldots,
P_k$ on $[n]$, and our computations amount to
finding the cardinality of the set $\mathcal{L}(P_i)$
of linear extensions of $P_i$.

The problem of computing linear extensions of general posets is
\#P-complete \cite{Brightwell1991}, so our task is an
intractable problem when $n$ grows large. However, 
for special classes of posets, and for moderate values of $n$,
the situation is not so bad. For example, in the  up-down 
analysis of Willbrand {\it et al.} (see Example \ref{ex.updwn}),
we need to count all permutations with a fixed descent set,
a task for which an explicit determinantal formula appears in Stanley 
\cite[page 69]{Stanley1997}. We refer to \cite{Brown2007}
for a detailed study of the combinatorics of these {\em up-down numbers}.

Likewise, there is an efficient (and easy-to-implement) method for the
computing quantities $\,|\, \tau^{-1} \bigl(\tau(\pi) \bigr)\,|\,$ for any graphical 
 graphical tubing test $\,\tau^*_{\mathcal{K}(G)}$, as defined in Section 5.
 Indeed, here the fan $\mathcal{F}$ is unimodular, and
 hence the posets $P_i$ are all trees. The special trees
 arising from a graph $G$ in this manner are known as {\em $G$-trees}
 \cite{Postnikov2005,Carr2004}. The $G$-tree of a permutation $\pi$
 is a representation of the poset $P_i$
 as a tree $\,T \,=\,\tau^*_{\mathcal{K}(G)}(\pi)\,$
 with the minimum value as the root and maximal values as the leaves. 
 Suppose the root of the tree $T$ has $k$
children, each of which is a root of a subtree $T^i$ for
$i=1,\ldots,k$. Writing $|T^i |$ for the number of nodes in
$T^i$, we have
$$ |\, \tau^{-1}(T ) \, |
\quad = \quad \binom{\sum_{i=1}^k |T^i |}{ |T^1|, \ldots, |
T^k|} \left( \prod_{i=1}^k |\tau^{-1}( T^{i})| \right).$$
This recursive formula translates into an efficient iterative
algorithm. Our implementation of this algorithm,
when $G$ is the $n$-cycle, is the workhorse
behind our computations in  \cite{cyclohedron}.
For a graph $G$, let $\nbhd(i)$ be the set of vertices $j$ such that there is an edge $(i,j)$ in $G$.

\begin{algorithm} \label{alg:TGMpermCount}(Permutation Counting)

\noindent {\em Input:} A data point $u$ as a descent permutation $\delta$ and a graph $G$.\\
\noindent {\em Output:} The number of permutations with the same signature as $\delta$, $|\, \tau^{-1}\tau(\pi(u)) \, |$.

\begin{verse}
{\bf Initialize:} \\
An indexed set of largest enclosing sets $LE_1=\dots=LE_n=\emptyset$, and counter $c=1$
\\
{\bf for} $\delta_i$ in $\delta$:\\
\quad Initialize $\ell$ an empty list of enclosed tree lengths\\
\quad $LE_{\delta_i}=\{\delta_i\}$\\
\quad {\bf for}  $j$ in $\nbhd(\delta_i)$:\\
\quad \quad {\bf if } $LE_j \neq \emptyset$ and $j \notin LE_{\delta_i}$:\\
\quad \quad \quad $LE_{\delta_i} = LE_{\delta_i} \disjointunion LE_j$\\
\quad \quad \quad append $|LE_j|$ to $\ell$\\
\quad  $c = c \cdot {\sum_i(\ell_i) \choose \ell}$\\
\quad {\bf for}  $j$ in $LE_{\delta_i}$:\\
\quad \quad $LE_j = LE_{\delta_i}$\\

{\bf Return} the permutation count $c$
\end{verse}
\end{algorithm}

In the remainder of this section we discuss our method
for performing these computations for an arbitrary convex rank test.
The test is specified (implicitly or explicitly)
 by a collection of posets  $P_1,\ldots,P_k$ on $[n]$.
From the given permutation, we identify the unique poset $P_i$
of which that permutation is a linear extension, and we construct
the {\em distributive lattice} $L(P_i)$ whose elements are the order ideals of
$P_i$. Recall that an {\em order ideal} of $P_i$ is  a subset $O$ of $[n]$
such that if $l \in O$ and $(k,l) \in P_i$ then $k \in O$. The set
of all order ideals is a distributive lattice with meet and join operations given
by set intersection $O \cap O'$ and set union $O \cup O'$. 

The distributive lattice $L(P_i)$ is a sublattice of the Boolean lattice $\,2^{[n]}$,
whose nodes are the $2^n$ subsets of $[n] = \{1,2,\ldots,n\}$,
and we represent $L(P_i)$ by its nodes and edges (cover relations) in $\,2^{[n]}$.
We write each edge in $2^{[n]}$ as a pair $(K,l)$
where $K \subset [n]$ and $l \in [n] \backslash K$. The edge 
in the Boolean lattice $2^{[n]}$ represented by the pair $(K,l)$
is the cover relation $\,K\, \subset \, K \cup \{l\}$.

Permutations in $S_n$ are in natural bijection with maximal chains
in the Boolean lattice $2^{[n]}$.  For example, the descent permutation $\delta=(4|2|1|3)$ corresponds to the maximal chain $\,\bigl(\emptyset, \{4\}, \{2, 4\}, \{1,2,4\}, \{1,2,3,4\}\bigr)\,$ in the Boolean lattice $\,2^{[4]}$. If the poset $P_i$ is the linear order $\delta$ then $L(P_i)$ is the subgraph of $2^{[4]}$
consisting of the five nodes in the chain and the four edges
$\,(\emptyset,4) ,\,(\{4\}, 2),\, (\{2,4\},1)\,$ and $\, (\{1,2,4\},3)\,$
which connect them. The maximal chains in $2^{[n]}$ that lie in the
sublattice $L(P_i)$ are precisely the permutations that
are linear extensions of $P_i$. Therefore our task
is to construct $L(P_i)$ and then count its maximal chains.

\begin{remark} \label{remlinex}
The linear extensions of the poset $P_i$ are in 
bijection with the maximal chains in the distributive lattice $L(P_i)$.
See \cite[Section 3.5]{Stanley1997} for further information on this bijection.
\end{remark}

In general, $L(P_i)$  is the graph whose nodes are those subsets of $[n]$
which are order ideals in $P_i$, and the edges
are $(K,l)$ where both $K$ and $K \cup \{l\} $ are order ideals in $P_i$.
Our strategy in computing the graph which represents $L(P_i)$
is as follows. We start with a given permutation $\delta$ which
lies in the class indexed by $P_i$. That permutation determines a 
maximal chain in $2^{[n]}$ which must lie in $L(P_i)$.
We then compute a certain closure of that subgraph in $2^{[n]}$
 with respect to the semigraphoid $\mathcal{M}$ under consideration.
This is precisely what is done in Algorithm 21 below.
Knowledge of the distributive lattice $L(P_i)$ solves our problem
since  the number of maximal chains of
$L(P_i)$ can be read easily from the representation of $L(P_i)$ in terms of nodes and edges.

\begin{algorithm} \label{DistLattice}(Building the Distributive Lattice)

\noindent {\em Input:} A data point as a descent permutation $\delta$ and a semigraphoid $\mathcal{M}$.\\
\noindent {\em Output:} A distributive lattice $L(P_i)$
 representing the class of $\delta$ in the convex rank test $\mathcal{M}$.

\begin{verse}
{\bf Initialize:} \\
A set of confirmed lattice nodes, $\,\mathbb{H} = 
\bigl\{ \emptyset, \{\delta_1\}, \{\delta_1, \delta_2\}, \dots, \{\delta_1, \dots, \delta_n\} \bigr\}$\\
A set of checked lattice edges, 
$\,E \,=\, \bigl\{ (\{\delta_1, \dots, \delta_{n-1}\}, \delta_n) \bigr\}$,\\
\ \ \ \ \ \ \ \ \ \ \ \ \ \ \ \ \ \ \ \ \ \ \ \ \ \ \ \ \ \ \ \ 
where each pair has the form \ (history, next position). \\
A stack of edges waiting to be checked: \\
 $\,W \,= \,
 \bigl[(\emptyset,\delta_1), (\{\delta_1\}, \delta_2), (\{\delta_1, \delta_2\},\delta_3), \dots,(\{\delta_1, \dots, \delta_{n-2}\}, \delta_{n-1}) \bigr]
$\\~\\
{\bf While} $W \neq \emptyset$:\\
\quad Pop $(H,i)$ from the stack $W$\\
\quad Add $(H,i)$ to $E$\\
\quad {\bf for}  $j$ such that $(H \union \{i\}, j) \in E$:\\
\quad \quad {\bf if} $i \indep j | H \in \mathcal{M}$:\\
\quad \quad \quad Add $(H,j)$ to $E$\\
\quad \quad \quad {\bf if} $H \union \{j\} \notin \mathbb{H}$:\\
\quad \quad \quad \quad Add $H \union \{j\}$ to $\mathbb{H}$\\
\quad \quad \quad \quad Push $(H \union \{j\}, i)$ onto $W$\\
{\bf Return} the distributive lattice $\,L(P_i)\,= \, \bigl(\mathbb{H},E \bigr) $\\
\end{verse}
\end{algorithm}

Our program for performing rank tests implements Algorithm \ref{DistLattice}.  It accepts a permutation $\delta$ and a rank test $\tau$, which may be specified either
\begin{itemize}
 \item by a list of posets $P_1,\ldots,P_k$ (pre-convex),
 \item or by a semigraphoid $\mathcal{M}$ (convex rank test),
 \item or by a submodular function $w : 2^{[n]} \rightarrow \R$,
 \item or by a collection $\mathcal{K}$ of subsets of $[n]$ (MSS),
 \item or by a graph $G$ on $[n]$ \  (graphical test).
\end{itemize}
The output of our program has two parts.
First, it gives the number $|\mathcal{L}(P_i)|$ of linear extensions,
where the poset $P_i$ represents the equivalence class of $S_n$ specified by the data $\pi$.
It also gives a representation of the distributive lattice $L(P_i)$, in a format
that can be read by the {\tt maple} package {\tt posets}  \cite{Stembridge2004}.
Our software for Algorithms \ref{alg:TGMpermCount} and \ref{DistLattice}
and, more generally, for
applying convex rank tests $\tau $ to data vectors $u \in \R^n$ is available at
$\, {\tt bio.math.berkeley.edu/ranktests/} $.

In closing let us give a concrete illustration of
our current ability to count linear extensions.
We computed the number of linear extensions 
of the Boolean poset $\,P  = 2^{[5]}\,$ consisting 
of all subsets of 
$\{1,2,3,4,5\}$. Our program ran in less than one second
on a laptop and found that 
$$ |L(2^{[5]})| \quad = \quad
14,807,804,035,657,359,360.$$
This computation was inspired by work 
in population genetics by 
Daniel Weinreich \cite{Weinreich2005} who reports the 
analogous calculation
for $P = 2^{[4]}$. 

\smallskip

\section*{Conclusions}

This work describes the connections among 
algebraic combinatorics, non-parametric statistics
and graphical models (statistical learning theory).
Specifically, we have proved the equivalence between semigraphoids and convex rank tests.  This result provides the background for the counterexamples given in \cite{counterexamples} and the rank tests which were applied to biological data in \cite{cyclohedron}.

\smallskip


\end{document}